\documentclass[12pt,reqno]{amsart}

\usepackage{pdfsync}

\usepackage{amssymb}

\usepackage{mathrsfs}

\DeclareMathAlphabet{\mathpzc}{OT1}{pzc}{m}{it}

\usepackage[all]{xy}

\entrymodifiers={+!!<0pt,\fontdimen22\textfont2>}

\def\cC{\mathscr{C}}

\def\cE{\mathscr{E}}
\def\cF{\mathscr{F}}

\def\cN{\mathscr{N}}

\def\cT{\mathscr{T}}

\def\BS{\mathbb{S}}

\def\fm{\mathfrak{m}}

\def\sD{\mathsf{D}}

\def\adots{\mathinner{\mkern1mu\raise1.0pt\vbox{\kern7.0pt\hbox{.}}\mkern2mu\raise4.0pt\hbox{.}\mkern2mu\raise7.0pt\hbox{.}\mkern1mu}}

\def\Cell{\operatorname{Cell}}

\def\CMreg{\operatorname{CMreg}}

\def\comp{\operatorname{comp}}
\def\D{\sD}

\def\Df{\D^{\operatorname{f}}}
\def\dim{\operatorname{dim}}

\def\Ext{\operatorname{Ext}}
\def\Extreg{\operatorname{Ext.\!reg}}

\def\Gr{\mathsf{Gr}}

\def\H{\operatorname{H}}

\def\Hom{\operatorname{Hom}}

\def\inf{\operatorname{inf}}

\newcommand\LTensor[1]{\overset{{\rm L}}{\underset{#1}{\otimes}}}

\def\opp{\operatorname{op}}

\def\RGammam{\operatorname{R}\!\Gamma_{\fm}}

\def\RHom{\operatorname{RHom}}

\def\sup{\operatorname{sup}}

\def\Tor{\operatorname{Tor}}
\def\tors{\operatorname{tors}}

\newtheorem{Lemma}{Lemma}[section]
\newtheorem{Theorem}[Lemma]{Theorem}
\newtheorem{Proposition}[Lemma]{Proposition}
\newtheorem{Corollary}[Lemma]{Corollary}

\theoremstyle{definition}
\newtheorem{Definition}[Lemma]{Definition}
\newtheorem{Setup}[Lemma]{Setup}

\newtheorem{Remark}[Lemma]{Remark}

\newtheorem{Example}[Lemma]{Example}

\newtheorem{Notation}[Lemma]{Notation}

\begin{document}

\setlength{\parindent}{0pt}
\setlength{\parskip}{7pt}

\title[Duality for cochain DG algebras]{Duality for cochain DG algebras}

\author{Peter J\o rgensen}
\address{School of Mathematics and Statistics,
Newcastle University, Newcastle upon Tyne NE1 7RU, United Kingdom}
\email{peter.jorgensen@ncl.ac.uk}
\urladdr{http://www.staff.ncl.ac.uk/peter.jorgensen}





\begin{abstract} 

This paper develops a duality theory for connected cochain DG
algebras, with particular emphasis on the non-com\-mu\-ta\-ti\-ve
aspects.

One of the main items is a dualizing DG module which induces a duality
between the derived categories of DG left-modules and DG right-modules
with finitely ge\-ne\-ra\-ted cohomology.

As an application, it is proved that if the canonical module $A /
A^{\geq 1}$ has a semi-free resolution where the cohomological degree
of the generators is bounded above, then the same is true for each DG
module with finitely generated cohomology.

\end{abstract}

\maketitle

\setcounter{section}{-1}
\section{Introduction}
\label{sec:introduction}

This paper develops a duality theory for connected cochain DG
algebras.  Some of the ingredients are dualizing DG modules, section
and completion functors, and local duality.

Particular emphasis is given to the non-commutative aspects of the
theory.  For instance, Theorem B below says that the dualizing DG
module defined in the paper induces a duality between the derived
categories of DG left-modules and DG right-modules with finitely
generated cohomology.

As an application, it is proved that if the canonical module $A /
A^{\geq 1}$ has a semi-free resolution where the cohomological degree
of the generators is bounded above, then the same is true for each DG
module with finitely generated cohomology.

\begin{Setup}
\label{set:A}
Throughout the paper, $A$ is a connected co\-cha\-in DG algebra over a
field $k$; that is, $A = A^{\geq 0}$ and $A^0 = k$.
\end{Setup}

See \cite[chp.\ 10]{BL} for an introduction to DG homological algebra.

Since $A$ is connected, it has a canonical DG bimodule $A / A^{\geq
1}$ which is also denoted by $k$.  The canonical module can be viewed
as a DG left-$A$-module ${}_{A}k$ or a DG right-$A$-module $k_A$, and
there are inclusions of lo\-ca\-li\-zing subcategories $\langle
{}_{A}k \rangle \hookrightarrow \sD(A)$ and $\langle k_A \rangle
\hookrightarrow \sD(A^{\opp})$ where $\sD(A)$ and $\sD(A^{\opp})$ are
the derived categories of DG left-$A$-modules, respectively DG
right-$A$-modules.  Under technical assumptions spelled out in Setup
\ref{set:2}, the inclusions have right-adjoint functors $\Gamma$ and
$\Gamma^{\opp}$ which behave like derived local section functors, and
the following is our first main result, see Theorem
\ref{thm:two-sided}.

\smallskip
\noindent
{\bf Theorem A. }
{\em
There is a single DG $A$-bimodule $F$ such that $\Gamma(-) = F
\LTensor{A} -$ and $\Gamma^{\opp}(-) = - \LTensor{A} F$.
}
\smallskip

The DG algebra $A$ may be very far from commutative, but the DG
bimodule $F$ behaves like a two-sided \v{C}ech complex for $A$ which
links DG left- and right-modules.  This is seen most clearly by
passing to the $k$-linear dual $D = \Hom_k(F,k)$ which will be called
a dualizing DG module of $A$.  The following is our second main
result, see Theorem \ref{thm:cor:two-sided}; like Theorem C below it
will be proved under the additional assumption that $\H(A)$ is
noetherian with a balanced dualizing complex.

\smallskip
\noindent
{\bf Theorem B. }
{\em
Let $\Df(A)$ and $\Df(A^{\opp})$ be the derived categories of DG
modules with finitely generated cohomology over $\H(A)$.  Then there
are quasi-inverse contravariant equivalences
\[
  \xymatrix{
    \Df(A) \ar[rrr]<1ex>^-{\RHom_A(-,D)}
    & & & \Df(A^{\opp}). \ar[lll]<1ex>^-{\RHom_{A^{\opp}}(-,D)}
           }
\]
}
\smallskip

As an application of the theory, we prove the following in Theorem
\ref{thm:regularities}.

\smallskip
\noindent
{\bf Theorem C. }
{\em
If ${}_{A}k$ has a semi-free resolution where the cohomological degree
of the generators is bounded above, then the same is true for each DG
left-$A$-module with finitely generated cohomology.
}
\smallskip

Note that despite the bound on the degree, there may be infinitely
many generators in each semi-free resolution of ${}_{A}k$.  For a
simple example, view $A = k[T]/(T^2)$ as a DG algebra with $T$ in
cohomological degree $1$ and $\partial = 0$.  Then the minimal
semi-free resolution of ${}_{A}k$ has all generators in degree $0$,
but there are infinitely many of them.  Hence each semi-free
resolution of ${}_{A}k$ has infinitely many generators and ${}_{A}k$
is not compact in the derived category.

The following describes three types of connected cochain DG algebras
where the results apply, see Section \ref{sec:examples}.
\begin{enumerate}

  \item  $\H(A)$ is noetherian AS regular.

\smallskip

  \item  $A$ is commutative in the DG sense and $\H(A)$ is noetherian.

\smallskip

  \item  $\dim_k \H(A) < \infty$.

\end{enumerate}
Between them, (i) and (ii) cover many DG algebras which arise in
practice, for instance as cochain DG algebras of topological spaces.
Note that in case (ii), Theorem A is trivial but Theorems B and C are
not.  In this case, the categories $\Df(A)$ and $\Df(A^{\opp})$ are
the same, but Theorem B says that this category has the non-trivial
property of being self dual.

Section \ref{sec:DG} summarises part of Dwyer and Greenlees's theory
of section and completion functors from \cite{DG}.  Section
\ref{sec:Shamir} considers the Greenlees spectral sequence \cite{G} in
a version given by Shamir \cite{Shamir}, evaluates it in the present
situation, and gives a technical consequence.  Section
\ref{sec:two-sided} proves Theorems A and B, Section
\ref{sec:regularities} proves Theorem C, and Section
\ref{sec:examples} provides some examples, not least by showing that
the theorems apply to the algebras described in (i)--(iii) above.

It should be mentioned that Mao and Wu \cite{MW} have provided some
technical tools which will be important in the proofs, and that $\Ext$
regularity (with the opposite sign) was studied in their paper under
the name width.  There is previous work on duality for DG algebras in
\cite{FIJ} and the more general $\BS$-algebras in \cite{DGI}.

\begin{Notation}
\label{not:blanket}
Opposite rings and DG algebras are denoted by the superscript
``$\opp$''.  Right (DG) modules are identified with left (DG) modules
over the opposite.  Subscripts are sometimes used to indicate left or
right (DG) module structures.

If $M$ is a DG module and $\ell$ an integer, then $M^{\geq \ell}$ is
the hard truncation with $(M^{\geq \ell})^n = 0$ for $n < \ell$.

Let $\sD$ denote the derived category of an abelian category or of
left DG modules over a DG algebra.  In the case of the standing DG
algebra $A$, let $\Df(A)$ be the full subcategory of objects $M \in
\sD(A)$ for which $\H(M)$ is finitely generated over $\H(A)$.

If $M$ is an object of a derived category $\sD$, then $\langle M
\rangle$ denotes the lo\-ca\-li\-zing subcategory generated by $M$.

If $\cN$ is a subcategory, then
\begin{align*}
  \cN^{\perp}   & = \{\, M \,|\, \Hom(N,M) = 0 \mbox{ for } N \in \cN \,\}, \\
  {}^{\perp}\!\cN & = \{\, M \,|\, \Hom(M,N) = 0 \mbox{ for } N \in \cN \,\}.
\end{align*}

The notation $(-)^*$ stands for $\Hom_A(-,A)$ or
$\Hom_{A^{\opp}}(-,A)$ and we write $(-)^{\vee} = \Hom_k(-,k)$.  These
functors interchange DG left- and right-$A$-modules.

For the theory of (balanced) dualizing complexes over connected graded
algebras see \cite{Y}.

If $M$ is a complex or a DG module, then we write
\[
  \inf M = \inf \{\, \ell \,|\, \H^\ell(M) \neq 0 \,\}, \;\;\;\;
  \sup M = \sup \{\, \ell \,|\, \H^\ell(M) \neq 0 \,\}.
\]
These are integers or $\pm \infty$.  Note that $\inf$ and $\sup$ of
the empty set are $\infty$ and $-\infty$, respectively, so $\inf(0) =
\infty$ and $\sup(0) = -\infty$.
\end{Notation}

\section{Dwyer-Greenlees theory}
\label{sec:DG}

\begin{Setup}
\label{set:blanket}
In this section and the next, we assume that $K$ is a K-projective DG
left-$A$-module which satisfies the following conditions as an object
of $\sD(A)$: It is compact and there is an equality of localizing
subcategories $\langle K \rangle = \langle {}_{A}k \rangle$.
\end{Setup}

\begin{Remark}
[Dwyer-Greenlees theory]
\label{rmk:DG}
The DG module $K$ can be used as an input for Dwyer and Greenlees's
theory from \cite{DG}.  Technically speaking, they only considered the
case of $K$ being a complex over a ring, but everything goes through
for a DG module over a DG algebra.  Let us give a brief recap of some
of their results.

Consider
\[
  \cE = \Hom_A(K,K)
\]
which is a DG algebra with multiplication given by composition.  Then
$K$ acquires the structure ${}_{A,\cE}K$ while $K^* = \RHom_A(K,A)$
has the structure $K^*_{A,\cE}$.  Define functors
\begin{align*}
  T(-) & = - \LTensor{\cE} K, \\
  E(-) & = \RHom_A(K,-) \simeq K^* \LTensor{A} -, \\
  C(-) & = \RHom_{\cE^{\opp}}(K^*,-)
\end{align*}
which form adjoint pairs $(T,E)$ and $(E,C)$ between $\sD(\cE^{\opp})$
and $\sD(A)$.

Set $\cN = \langle {}_{A}k \rangle^{\perp} = \langle {}_{A}K
\rangle^{\perp}$ in $\sD(A)$; in terms of these null modules we define
the torsion and the complete DG modules by
\[
  \sD^{\tors}(A) = {}^{\perp}\!\cN, \;\;\; \sD^{\comp}(A) = \cN^{\perp}.
\]
Note that
\[
  \sD^{\tors}(A) = \langle {}_{A}k \rangle = \langle {}_{A}K \rangle.
\]
There are pairs of quasi-inverse equivalences of categories as
follows.
\[
  \xymatrix{
    \sD^{\comp}(A) \ar[rr]<1ex>^-{E}
    & & \sD(\cE^{\opp}) \ar[ll]<1ex>^-{C} \ar[rr]<1ex>^-{T}
    & & \sD^{\tors}(A) \ar[ll]<1ex>^-{E}
           }
\]
In particular, $EC$ and $ET$ are equivalent to the identity functor on
$\sD(\cE^{\opp})$, so if we set
\[
  \Gamma = TE, \;\;\; \Lambda = CE
\]
then we get endofunctors of $\sD(A)$ which form an adjoint pair
$(\Gamma,\Lambda)$ and satisfy
\[
  \Gamma^2 \simeq \Gamma, \;\;\;
  \Lambda^2 \simeq \Lambda, \;\;\;
  \Gamma\Lambda \simeq \Gamma, \;\;\;
  \Lambda\Gamma \simeq \Lambda.
\]
These functors are adjoints of inclusions as follows, where
left-adjoints are displayed above right-adjoints.
\[
  \xymatrix{
    \sD^{\comp}(A) \ar[rr]<-1ex>_-{\operatorname{inc}}
    & & \sD(A) \ar[ll]<-1ex>_-{\Lambda} \ar[rr]<-1ex>_-{\Gamma}
    & & \sD^{\tors}(A) \ar[ll]<-1ex>_-{\operatorname{inc}}
           }
\]
Note that the counit and unit, $\Gamma(-)
\stackrel{\epsilon}{\longrightarrow} (-)$ and $(-)
\stackrel{\eta}{\longrightarrow} \Lambda(-)$, are
$K$-equivalences, that is, they become isomorphisms when the functor
$\Hom_{\sD(A)}(\Sigma^\ell K,-)$ is applied.  Equivalently, their
mapping cones are in $\langle {}_{A}k \rangle^{\perp}$.  Along with
$\Gamma(-) \in \sD^{\tors}(A)$ and $\Lambda(-) \in \sD^{\comp}(A)$,
this characterizes them up to unique isomorphism.

It is useful to remark that in particular, for $M \in
\sD^{\tors}(A) = \langle {}_{A}k \rangle$, the counit morphism
$\Gamma(M) \stackrel{\epsilon_M}{\longrightarrow} M$ is an isomorphism,
and for $M \in \sD^{\comp}(A)$, the unit morphism $M
\stackrel{\eta_M}{\longrightarrow} \Lambda(M)$ is an isomorphism.  For $M
\in \langle {}_{A}k \rangle^{\perp}$, we get $\Gamma(M) = 0$.
\end{Remark}

\begin{Definition}
\label{def:F}
We will write
\[
  F = K^* \LTensor{\cE} K, \;\;\;
  D = F^{\vee}
\]
and refer to $D$ as a dualizing DG module of $A$.

In a more laborious notation we have $F = K^*_{A,\cE}
\LTensor{\cE} {}_{A,\cE}K$, so $F$ has the structure ${}_{A}F_{A}$ and
$D$ the structure ${}_{A}D_{A}$.  It is easy to check that
\[
  \Gamma(-)  = F \LTensor{A} -, \;\;\;
  \Lambda(-) = \RHom_A(F,-)
\]
and adjointness yields
\begin{equation}
\label{equ:local_duality}
  \Gamma(-)^{\vee} = \RHom_A(-,D).
\end{equation}
The DG module $F$ plays the role of the \v{C}ech complex and $\Gamma$
and $\Lambda$ behave like derived local section and completion
functors.  Equation \eqref{equ:local_duality} is the local duality
formula.
\end{Definition}

\begin{Remark}
\label{rmk:two-sided}
Since $\Gamma$ and $\Lambda$ are given as derived $\otimes$ and $\Hom$
with the DG $A$-bimodule $F$, they can be applied to DG $A$-bimodules
and this will give new DG $A$-bimodules.

Specifically, $\Gamma({}_{A}M_{A}) = {}_{A}F_{A} \LTensor{A}
{}_{A}M_{A}$ has a left-structure which comes from the left-structure
of $F$ and a right-structure which comes from the right-structure of
$M$.  And $\Lambda({}_{A}M_{A}) = \RHom_A({}_{A}F_{A} , {}_{A}M_{A})$
has a left-structure which comes from the right-structure of $F$ and a
right-structure which comes from the right-structure of $M$.

It is easy to check that when the functors are applied to DG
$A$-bimodules, the counit and unit, $\Gamma(-)
\stackrel{\epsilon}{\longrightarrow} (-)$ and $(-)
\stackrel{\eta}{\longrightarrow} \Lambda(-)$, can be viewed as
morphisms in $\sD(A^e)$, the derived category of DG $A$-bimodules.
\end{Remark}

\section{The Greenlees spectral sequence in a version given by Shamir}
\label{sec:Shamir}

\begin{Remark}
\label{rmk:spectral_sequence}
In this remark, assume that $\H(A)$ is noetherian.

The Greenlees spectral sequence was originally given for group
cohomology in \cite[thm.\ 2.1]{G}.  It was developed further by
Benson, Dwyer, Greenlees, Iyengar, and Shamir in \cite{BG},
\cite{DGI}, and \cite{Shamir}.  The most general version is given by
Shamir in \cite{Shamir}; we will apply it to the situation at hand.


The cohomology $\H(A)$ is a connected graded $k$-algebra with graded
maximal ideal $\fm = \H^{\geq 1}(A)$.  Let $\cT$ denote the
$\fm$-torsion graded left-$\H(A)$-modules, that is, the graded
mo\-du\-les such that each element $t$ has $\fm^{\ell}t = 0$ for $\ell
\gg 0$.  Then $\cT$ is a hereditary torsion class in the abelian
category $\Gr \H(A)$ of graded left-$\H(A)$-modules, in the sense of
\cite[def.\ 3.1]{Shamir}.  For $X \in \Gr \H(A)$, view $X$ as an
object of the derived category $\sD(\Gr \H(A))$ and, using the
notation of \cite[def.\ 2.1]{Shamir}, consider a morphism
$\Cell^{\H(A)}_{\cT}(X) \stackrel{\eta}{\longrightarrow} X$ in
$\sD(\Gr \H(A))$ characterized by the properties that
$\Cell^{\H(A)}_{\cT}(X) \in \langle \cT \rangle$ and $\Hom_{\sD(\Gr
\H(A))}(\Sigma^{\ell} T,\eta)$ is an isomorphism for each integer
$\ell$ and $T \in \cT$.  These properties determine $\eta$ up to
unique isomorphism, and using that $\H(A)$ is noetherian, it is not
hard to show that $\eta$ can be obtained as the canonical morphism
$\RGammam X \rightarrow X$ where $\RGammam$ is the functor on the
derived category which underlies local cohomology; see \cite[sec.\
4]{Y}.

Consider the class
\[
  \cC = \{\, C \in \sD(A) \,|\, \H(C) \in \cT \,\}.
\]
Shamir refers to objects of $\cC$ as $\cT$-cellular and objects of
$\cC^{\perp}$ as $\cT$-null; see for instance \cite[p.\ 1 and defs.\
2.1 and 2.3]{Shamir}.  For each DG left-$A$-module $M$, Shamir obtains
in \cite[lem.\ 5.4]{Shamir} a distinguished triangle $C
\rightarrow M \rightarrow N$ in $\sD(A)$ with $N \in \cC^{\perp}$
and a spectral sequence
\[
  E^2_{p,q} = \H_{p,q}(\RGammam(\H\!M)) \Rightarrow \H_{p+q}(C).
\]
On the left hand side, $p$ comes from the numbering of the modules
in an exact couple and $q$ is an internal degree; see \cite[proof of
lem.\ 5.4]{Shamir}.
A consequence is that $p$ is homological degree
along the complex $\RGammam$ and $q$ is graded degree along the graded
module $\H\!M$.  The sequence can hence also be written
\begin{equation}
\label{equ:spectral_sequence}
  E^2_{p,q} = \H_{\fm}^{-p}(\H\!M)_q \Rightarrow \H_{p+q}(C)
\end{equation}
where $\H_{\fm}^\ell = \H^\ell \circ \RGammam$ is local cohomology;
see \cite[sec.\ 4]{Y}.

The spectral sequence is conditionally convergent to the colimit;
compare \cite[last part of proof of lem.\ 5.4]{Shamir} with
\cite[def.\ 5.10]{Boardman}.  In fact, the $p$ in \cite{Shamir}
corresponds to the $s$ in \cite[eq.\ (0.1)]{Boardman}, except that
they have opposite signs.  Now note that $E^2_{p,*} = 0$ for $p > 0$
by construction, so the spectral sequence is a half-plane spectral
sequence in the sense of \cite[sec.\ 7]{Boardman}.  By \cite[thm.\
7.1]{Boardman}, to get strong convergence, all we need is to check
$\operatorname{R}\!E_{\infty} = 0$ in the notation of \cite{Boardman}.

In particular, using \cite[rmk.\ after thm.\ 7.1]{Boardman}, the
spectral sequence \eqref{equ:spectral_sequence} is strongly convergent
if
\begin{equation}
\label{equ:finite_dimensions}
  \dim_k E^2_{p,q} = \dim_k \H_{\fm}^{-p}(\H\!M)_q < \infty
\end{equation}
for all $p,q$.
\end{Remark}

\begin{Lemma}
\label{lem:bounded_above}
If $M \in \sD(A)$ has $\H^\ell(M) = 0$ for $\ell \gg 0$ then $M \in
\langle {}_{A}k \rangle$.
\end{Lemma}

\begin{proof}
Using \cite[sec.\ 1.5]{MW} to truncate, we can suppose $M^\ell = 0$
for $\ell \gg 0$ and desuspending if necessary, we can suppose $M^\ell
= 0$ for $\ell > 0$ so $M^{\geq 1} = 0$.  There is a direct system
\[
  M^{\geq 0}
  \rightarrow M^{\geq -1}
  \rightarrow M^{\geq -2}
  \rightarrow \cdots
\]
in $\sD(A)$ with homotopy colimit $M$, so it is enough to see $M^{\geq
n} \in \langle {}_{A}k \rangle$ for each $n$.  However, there are
distinguished triangles
\[
  M^{\geq n+1} \rightarrow M^{\geq n} \rightarrow N(n)
\]
where $\H(N(n))$ is concentrated in cohomological degree $n$, and it
is easy to check that hence $N(n) \cong \coprod \Sigma^{-n}k$ so $N(n)
\in \langle {}_{A}k \rangle$.  Induction starting with $M^{\geq 1} =
0$ gives $M^{\geq n} \in \langle {}_{A}k \rangle$ for each $n$ as
desired.
\end{proof}

The following proof uses the methods of \cite[sec.\ 6]{Shamir}. 

\begin{Theorem}
\label{thm:spectral_sequence}
Assume (in addition to Setup \ref{set:blanket}) that $\H(A)$ is
noetherian with a balanced dualizing complex.  For $M \in \Df(A)$
there is a spectral sequence which is strongly convergent in the sense
of \cite[def.\ 5.2]{Boardman}, 
\[
  E^2_{p,q} = \H_{\fm}^{-p}(\H\!M)_q \Rightarrow \H_{p+q}(\Gamma M).
\]
\end{Theorem}

\begin{proof}
The condition in Equation \eqref{equ:finite_dimensions} holds because
$\H(A)$ is noetherian with a balanced dualizing complex and $\H(M)$ is
finitely generated; combine \cite[thms.\ 5.1 and 6.3]{VdB}.  Hence
the spectral sequence (\ref{equ:spectral_sequence}) is strongly
convergent by the observation at the end of Remark
\ref{rmk:spectral_sequence}. 

To complete the proof, we must see $C \cong \Gamma M$ where $C$ is the
object in (\ref{equ:spectral_sequence}).  There is a distinguished
triangle $C \rightarrow M \rightarrow N$, so by \cite[prop.\
2.4]{Shamir} it is enough to see $C \in \langle {}_{A}k \rangle$ and
$N \in \langle {}_{A}k \rangle^{\perp}$.

For the latter, it suffices to see $\Hom_A(\Sigma^\ell k,N) = 0$ for
each $\ell$, and this is clear because $N \in \cC^{\perp}$ in the
notation from Remark \ref{rmk:spectral_sequence}.

For the former, note that $\H(A)$ is noetherian with a balanced
dualizing complex and $\H(M)$ is finitely generated.  Hence
$\H_{\fm}^{>n}(\H\!M) = 0$ for some $n$, and for each $p$ the graded
module $\H_{\fm}^p(\H\!M)$ is zero in sufficiently high degree; this
is by \cite[thms.\ 5.1 and 6.3]{VdB} again.  The degree in question
stems from the cohomological grading of $\H(M)$ so we learn
$\H_{\fm}^p(\H\!M)_q = 0$ for $q \ll 0$ since $q$ figures as a
subscript, hence with a sign change.  So $E^2_{p,q}$ is concentrated
in a vertical strip which is bounded below, and strong convergence to
$\H_{p+q}(C)$ implies $\H_\ell(C) = 0$ for $\ell \ll 0$.  But then $C
\in \langle {}_{A}k \rangle$ by Lemma \ref{lem:bounded_above}.
\end{proof}

\begin{Remark}
Note that it is easy to show that $C \in \langle {}_{A}k \rangle$
implies $\H(C) \in \cT$.  If we also had the opposite implication,
then we could conclude that $\Gamma(M)$ was $\Cell_{\cT}^A(M)$ in the
notation of \cite{Shamir}, and obtain Theorem
\ref{thm:spectral_sequence} as a special case of \cite[thm.\
1]{Shamir}.
\end{Remark}

\begin{Corollary}
\label{cor:spectral_sequence}
Assume that $\H(A)$ is noetherian with a balanced du\-a\-li\-zing
complex.  Then $M \in \Df(A)$ implies $(\Gamma M)^{\vee} \in
\Df(A^{\opp})$.
\end{Corollary}

\begin{proof}
We must show that if $\H(M)$ is finitely generated over $\H(A)$ then
$\H(\Gamma M)^{\vee}$ is finitely generated over $\H(A)^{\opp}$.

As indicated in Remark \ref{rmk:spectral_sequence}, in the spectral
sequence, $q$ is internal degree.  The same hence applies to the
spectral sequence of Theorem \ref{thm:spectral_sequence}.  But by
\cite[thms.\ 5.1 and 6.3]{VdB}, the graded modules
$\H_{\fm}^{-p}(\H\!M)$ are the $k$-linear duals of finitely generated
$\H(A)^{\opp}$-modules, and only finitely many of them are non-zero.

So the terms $E^2_{p,*}$ are the $k$-linear duals of finitely many
finitely ge\-ne\-ra\-ted $\H(A)^{\opp}$-modules whence the same is
true for the terms $E^{\infty}_{p,*}$.  Hence $\H(\Gamma M)$ has a
filtration where the quotients are the $k$-linear duals of finitely
many finitely generated $\H(A)^{\opp}$-modules.  This proves the
result.
\end{proof}

\section{Properties of the \v{C}ech and dualizing DG modules}
\label{sec:two-sided}

\begin{Setup}
\label{set:2}
In this section and the next, we assume that ${}_{A}K$ and $L_A$ are
K-projective DG $A$-modules which satisfy the following conditions as
objects of $\sD(A)$ and $\sD(A^{\opp})$:
\begin{itemize}

  \item   ${}_{A}K$ is compact, $\langle {}_{A}K \rangle = \langle
{}_{A}k \rangle$, and $K^*_A \in \langle k_A \rangle$,

\smallskip

  \item  $L_A$ is compact, $\langle L_A \rangle = \langle
k_A \rangle$, and ${}_{A}L^* \in \langle {}_{A}k \rangle$.

\end{itemize}
\end{Setup}

\begin{Remark}
\label{rmk:L}
The DG module $K$ can be used as input for the theory of the previous
sections.  In particular, Section \ref{sec:DG} used $K$ to define
various objects which will be important: $\cE$, $F$, $D$, $\Gamma$,
$\Lambda$.

Similarly, $L$ can be used as input for the theory applied to
$A^{\opp}$, that is, to DG right-$A$-modules.  In this case we get the
endomorphism DG algebra
\[
  \cF = \Hom_{A^{\opp}}(L,L),
\]
the DG module $L$ acquires the structure ${}_{\cF}L_A$, and we can
define the DG $A$-bimodules
\[
  G = L^* \LTensor{\cF} L, \;\;\; E = G^{\vee}
\]
along with the functors
\[
  \Gamma^{\opp}(-)  = - \LTensor{A} G, \;\;\;
  \Lambda^{\opp}(-) = \RHom_{A^{\opp}}(G,-).
\]
Then $\cF$, $G$, $E$, $\Gamma^{\opp}$, $\Lambda^{\opp}$ are the right
handed versions of $\cE$, $F$, $D$, $\Gamma$, $\Lambda$.
\end{Remark}

The following is Theorems A and B of the introduction.

\begin{Theorem}
\label{thm:two-sided}
We have $F \cong G$ and $D \cong E$ in the derived category $\sD(A^e)$
of DG $A$-bimodules. 
\end{Theorem}

\begin{proof}
We know that ${}_{\cF}L$ is built from ${}_{\cF}\cF$ using
(de)suspensions, distinguished triangles, coproducts, and direct
summands.  The functor ${}_{A}L_{\cF}^* \LTensor{\cF} -$ preserves
these operations and $\langle {}_{A}k \rangle$ is closed under them,
so ${}_{A}L^* \in \langle {}_{A}k \rangle$ implies ${}_{A}G =
{}_{A}L^*_{\cF} \LTensor{\cF} {}_{\cF}L \in \langle {}_{A}k
\rangle$.  By symmetry, $F_A \in \langle k_A \rangle$.

However, we have
\[
  {}_{A}G_{A}
  \stackrel{\epsilon_G}{\longleftarrow} \Gamma({}_{A}G_{A})
  = {}_{A}F_{A} \LTensor{A} {}_{A}G_{A}
  = \Gamma^{\opp}({}_{A}F_{A})
  \stackrel{\epsilon^{\opp}_F}{\longrightarrow} {}_{A}F_{A}
\]
where the counit morphisms $\epsilon_G$ and $\epsilon^{\opp}_F$ are
morphisms in $\sD(A^e)$ as explained in Remark \ref{rmk:two-sided}.
Now, ${}_{A}G$ is in $\langle {}_{A}k \rangle$ so by the last
paragraph of Remark \ref{rmk:DG}, if we forget the
right-$A$-structures, then $\epsilon_G$ is an isomorphism in $\sD(A)$.
This just means that its cohomology is bijective whence $\epsilon_G$
itself is an isomorphism in $\sD(A^e)$.  By symmetry,
$\epsilon^{\opp}_F$ is an isomorphism in $\sD(A^e)$ and the
proposition follows.
\end{proof}

\begin{Theorem}
\label{thm:cor:two-sided}
Assume (in addition to Setup \ref{set:2}) that $\H(A)$ is
noetherian with a balanced dualizing complex.  Then there are
quasi-inverse contravariant equivalences 
\[
  \xymatrix{
    \Df(A) \ar[rrr]<1ex>^-{\RHom_A(-,D)}
    & & & \Df(A^{\opp}). \ar[lll]<1ex>^-{\RHom_{A^{\opp}}(-,D)}
           }
\]
\end{Theorem}

\begin{proof}
Definition \ref{def:F}, Remark \ref{rmk:L}, and Theorem
\ref{thm:two-sided} show that the two functors in the theorem are
$\Gamma(-)^{\vee}$ and $\Gamma^{\opp}(-)^{\vee}$.  They take values in
the correct categories by Corollary \ref{cor:spectral_sequence} and
its analogue for $\Gamma^{\opp}$.

To see that the functors are quasi-inverse equivalences, first observe
that by adjointness,
\[
  \Gamma(-)^{\vee}
  = (F \LTensor{A} -)^{\vee}
  \simeq \RHom_{A^{\opp}}(F,(-)^{\vee})
  = \Lambda^{\opp}((-)^{\vee}).
\]
This gives the first of the following natural isomorphisms for $M \in 
\Df(A)$.
\[
  \Gamma^{\opp}(\Gamma(M)^{\vee})^{\vee}
  \cong \Gamma^{\opp}(\Lambda^{\opp}(M^{\vee}))^{\vee}
  \stackrel{\rm (a)}{\cong} \Gamma^{\opp}(M^{\vee})^{\vee}
  \stackrel{\rm (b)}{\cong} M^{\vee \vee}
  \cong M.
\]
Here (a) is because $\Gamma^{\opp}\Lambda^{\opp} \simeq
\Gamma^{\opp}$ by Remark \ref{rmk:DG}, and (b) is because when
$\H(M)$ is finitely generated, it has $\H^{\ell}(M) = 0$ for $\ell \ll
0$ whence $\H^{\ell}(M^{\vee}) = 0$ for $\ell \gg 0$; hence $M^{\vee}
\in \langle k_A \rangle$ by the right-module version of Lemma
\ref{lem:bounded_above} and so $\Gamma^{\opp}(M^{\vee}) \cong
M^{\vee}$.

The reverse composition of functors is handled by symmetry.
\end{proof}

\section{An application to Ext regularity}
\label{sec:regularities}

\begin{Definition}
For $M \in \sD(A)$ we define the $\Ext$ and
Ca\-stel\-nu\-o\-vo-Mum\-ford regularities by
\[
  \Extreg M  = - \inf \RHom_A(M,k),  \;\;\;
  \CMreg M   = \sup \Gamma M,
\]
and similarly for $M \in \sD(A^{\opp})$.

Note that $\Extreg(0) = \CMreg(0) = -\infty$; see the last part of
Notation \ref{not:blanket}.
\end{Definition}

\begin{Remark}
\label{rmk:regularities}
If $M \in \sD(A)$ has $\H^\ell(M) = 0$ for $\ell \ll 0$, then it
follows from \cite[prop.\ 2.4]{MW} that $M$ has a minimal semi-free
resolution $P$ with generators between cohomological degrees $\inf M$
and $\Extreg M$.  That is, if we write $i = \inf M$ and $r = \Extreg
M$, then
\begin{equation}
\label{equ:P}
  P^{\natural} = \coprod_{-r \leq \ell \leq -i} \Sigma^{\ell}(A^{\natural})^{(\beta_{\ell})}
\end{equation}
where $\natural$ sends DG modules to graded modules by forgetting the
dif\-fe\-ren\-ti\-al and $(\beta_{\ell})$ indicates a coproduct.
If, additionally, $\H(M) \neq 0$, then $\inf M$ is a finite number and
$P$ has at least one generator.  Hence
\begin{align}
\nonumber
  & \mbox{$\H^{\ell}(M) = 0$ for $\ell \ll 0$ and $\H(M) \neq 0$} \\
\label{equ:inf_and_Extreg}
  & \hspace{21ex} \Rightarrow \; -\infty < \inf M \leq \Extreg M.
\end{align}

If $\H(A)$ is noetherian with a balanced dualizing complex, then
Equation \eqref{equ:local_duality} in Definition \ref{def:F} along with
Theorem \ref{thm:cor:two-sided} give
\begin{equation}
\label{equ:CMreg_finite}
  M \in \Df(A) \mbox{ and } \H(M) \neq 0
  \; \Rightarrow \; -\infty < \CMreg M < \infty.
\end{equation}

By considering $k_A \LTensor{A} {}_{A}k$, one proves $\Extreg(k_A) =
\Extreg({}_{A}k)$, and this common number will be denoted by $\Extreg
k$.  Equation \eqref{equ:inf_and_Extreg} implies
\[
  0 \leq \Extreg k.
\]
By using Theorem \ref{thm:two-sided} we get $\Gamma^{\opp}(A)
= A \LTensor{A} F \cong F_A$ and $\Gamma(A) = F \LTensor{A} A \cong
{}_{A}F$, so $\CMreg(A_A) = \CMreg({}_{A}A)$, and this common number will
be denoted by $\CMreg A$.
\end{Remark}

\begin{Definition}
[He and Wu {\cite[def.\ 2.1]{HW2}}]
A DG $A$-module $M$ is Koszul if it has a semi-free resolution $P$
all of whose basis elements are in degree $0$.

The DG algebra $A$ is Koszul if ${}_{A}k$ is a Koszul DG module.
\end{Definition}

\begin{Remark}
\label{rmk:Koszul}
If $M$ is a DG $A$-module with $\H^\ell(M) = 0$ for $\ell \ll 0$, then
it is immediate from Remark \ref{rmk:regularities} that it is Koszul
precisely if $\H(M) = 0$ or $\inf M = \Extreg M = 0$.

Consequently, the DG algebra $A$ is Koszul precisely if $\Extreg k =
0$. 
\end{Remark}

\begin{Lemma}
\label{lem:formula}
Suppose that $M \in \sD(A)$ has $\H^\ell(M) = 0$ for $\ell \ll 0$ and
$\dim_k \H^\ell(M) < \infty$ for each $\ell$.  Then $\Lambda(M) \cong M$.
\end{Lemma}

\begin{proof}
The assumptions on $M$ imply $M \cong \RHom_k(M^{\vee},k)$.  This
gives the first of the following isomorphisms, and the second one is
by adjointness. 
\begin{align*}
  \Lambda(M)
  & = \RHom_A(F,M) \\
  & \cong \RHom_A(F,\RHom_k(M^{\vee},k)) \\
  & \cong \RHom_k(M^{\vee} \LTensor{A} F,k) \\
  & = (M^{\vee} \LTensor{A} F)^{\vee} \\
  & = \Gamma^{\opp}(M^{\vee})^{\vee} \\
  & \stackrel{\rm (a)}{\cong} (M^{\vee})^{\vee} \\
  & \cong M
\end{align*}
Here (a) follows from the right-module version of Lemma
\ref{lem:bounded_above} and the final paragraph of Remark
\ref{rmk:DG}. 

Note that the proof uses the two-sided theory of Section
\ref{sec:two-sided}: It is necessary to know that $\Lambda$ and
$\Gamma^{\opp}$ are given by formulae involving the {\em same} DG
bimodule $F$.
\end{proof}

The following is a DG version of \cite[thms.\ 2.5 and 2.6]{J}.

\begin{Proposition}
\label{pro:regularities}
Let $M \in \sD(A)$ have $\H^\ell(M) = 0$ for $\ell \ll 0$ and $\H(M)
\neq 0$.  Then
\begin{enumerate}

  \item  $\CMreg M \neq -\infty$.

\smallskip

  \item  $\Extreg M \leq \CMreg M + \Extreg k$.

\smallskip

  \item  $\CMreg M \leq \Extreg M + \CMreg A$.

\end{enumerate}
\end{Proposition}

\begin{proof}
(i) Observe that $\Lambda k \cong k$ by Lemma \ref{lem:formula}, so
\begin{align}
\nonumber
  \RHom_A(M,k)
  & \cong \RHom_A(M , \Lambda k) \\
\nonumber
  & = \RHom_A(M , \RHom_A(F,k)) \\
\nonumber
  & \stackrel{\rm (a)}{\cong} \RHom_A(F \LTensor{A} M , k) \\
\label{equ:formula}
  & = \RHom_A(\Gamma M,k)
\end{align}
where (a) is by adjointness.  Hence
\begin{equation}
\label{equ:z}
  \Extreg M = -\inf \RHom_A(\Gamma M,k).
\end{equation}
Now, $\CMreg M = \sup \Gamma M = -\infty$ would mean $\Gamma M = 0$.
By Equation \eqref{equ:z} this would imply $\Extreg M = -\infty$, but
this is false by Equation \eqref{equ:inf_and_Extreg} in Remark
\ref{rmk:regularities}. 

\smallskip
(ii) By part (i) and Equation \eqref{equ:inf_and_Extreg} we have
$\CMreg M$ and $\Extreg k$ different from $-\infty$.  Hence, despite
the potential for either regularity to be $\infty$, the right hand
side of the inequality in the proposition makes sense because it does
not read $\infty - \infty$.

Set $X = \Gamma M$ and let $P$ be a minimal semi-free resolution of
$k_A$.  Then
\[
  \Extreg M
  \stackrel{\rm (b)}{=} -\inf \RHom_A(X,k) 
  \stackrel{\rm (c)}{=} -\inf \Hom_A(X,P^{\vee}) 
  = (*)
\]
where (b) is Equation \eqref{equ:z} and (c) is because $P^{\vee}$ is a
K-injective resolution of ${}_{A}k$.

We have $\sup X = \CMreg M$ so by truncation we can suppose $X^\ell =
0$ for $\ell > \CMreg M$, cf.\ \cite[1.6]{MW}, and so
\[
  \mbox{ $(X^{\vee})^j = 0$ for $j < -\CMreg M$. }
\]
Write $i = \inf k = 0$ and $r = \Extreg k$.  Then $P$ satisfies
Equation \eqref{equ:P} in Remark \ref{rmk:regularities} and a
computation shows 
\[
  \Hom_A(X,P^{\vee})^{\natural}
  \cong \prod_{-r \leq \ell} \Sigma^{-\ell}((X^{\vee})^{\natural})^{\beta_{\ell}}
\]
where the power $\beta_{\ell}$ indicates a product.  The last two
equations imply
\[
  (*) \leq \CMreg M + r = \CMreg M + \Extreg k
\]
as desired.

\smallskip
(iii)  Note that the right hand side of the inequality makes sense
again, for the same reason as in part (ii).

Let $P$ be a minimal semi-free resolution of $M$.  Then 
\[
  \CMreg M
  = \sup \Gamma M
  = \sup F \LTensor{A} M
  = \sup F \underset{A}{\otimes} P
  = (**).
\]

As noted in Remark \ref{rmk:regularities} we have $F_A \cong
\Gamma^{\opp}(A)$ and since $\CMreg A = \sup \Gamma^{\opp}(A)$
we can suppose by truncation that
\[
  \mbox{ $F^j = 0$ for $j > \CMreg A$. }
\]
Write $i = \inf M$ and $r = \Extreg M$.  Then $P$ satisfies Equation
\eqref{equ:P} in Remark \ref{rmk:regularities} and a computation shows
\[
  (F \underset{A}{\otimes} P)^{\natural}
  \cong \coprod_{-r \leq n} \Sigma^n(F^{\natural})^{(\beta_n)}.
\]
The last two equations imply
\[
  (**) \leq r + \CMreg A = \Extreg M + \CMreg A
\]
as desired.
\end{proof}

Part (i) of the following establishes Theorem C of the introduction
while (ii) is a DG version of \cite[cor.\ 2.9]{J}.

\begin{Theorem}
\label{thm:regularities}
Assume (in addition to Setup \ref{set:2}) that $\H(A)$ is
noe\-the\-ri\-an with a balanced dualizing complex.  Let $M \in
\Df(A)$ have $\H(M) \neq 0$. 
\begin{enumerate}

  \item  If $\Extreg k < \infty$ then $\Extreg M < \infty$.

\smallskip

  \item  If $A$ is a Koszul DG algebra and $\CMreg M \leq t$ for an
integer $t$, then $\Sigma^t(M^{\geq t})$ is a Koszul DG module.

\end{enumerate}
\end{Theorem}

\begin{proof}
(i) follows by combining Equation \eqref{equ:CMreg_finite} in Remark
\ref{rmk:regularities} with Proposition \ref{pro:regularities}(ii).

As for (ii), it holds trivially if $\H(\Sigma^t(M^{\geq t})) = 0$,
so suppose that we have $\H(\Sigma^t(M^{\geq t})) \neq 0$.

There is a short exact sequence of DG modules $0 \rightarrow M^{\geq
  t} \rightarrow M \rightarrow M/M^{\geq t} \rightarrow 0$ which
induces a distinguished triangle $\Sigma^{-1}(M/M^{\geq t})
\rightarrow M^{\geq t} \rightarrow M$ in $\sD(A)$, and hence a
distinguished triangle
\[
  \Gamma(\Sigma^{-1}(M/M^{\geq t}))
    \rightarrow \Gamma(M^{\geq t})
    \rightarrow \Gamma M
\]
in $\sD(A)$.

Lemma \ref{lem:bounded_above} and the last paragraph of Remark
\ref{rmk:DG} imply the isomorphism $\Gamma(\Sigma^{-1}(M/M^{\geq t}))
\cong \Sigma^{-1}(M/M^{\geq t})$, so $\sup
\Gamma(\Sigma^{-1}(M/M^{\geq t})) \leq t$.  By assumption, $\sup
\Gamma M = \CMreg M \leq t$.  So the distinguished triangle implies
$\sup \Gamma(M^{\geq t}) \leq t$.  Hence $\sup \Gamma(\Sigma^t(M^{\geq
  t})) \leq 0$, that is $\CMreg \Sigma^t(M^{\geq t}) \leq 0$.  But
then $\Extreg \Sigma^t(M^{\geq t}) \leq 0$ by Proposition
\ref{pro:regularities}(ii).

On the other hand, it is clear that $\inf \Sigma^t(M^{\geq t}) \geq
0$, and Equation \eqref{equ:inf_and_Extreg} in Remark
\ref{rmk:regularities} now implies $\inf \Sigma^t(M^{\geq t}) =
\Extreg \Sigma^t(M^{\geq t}) = 0$, so $\Sigma^t(M^{\geq t})$ is a
Koszul DG module.
\end{proof}

\section{Examples}
\label{sec:examples}

Recall from Setup \ref{set:A} that $A$ is a connected cochain DG
algebra over a field $k$.

\begin{Example}
\label{exa:3}
If $\H(A)$ is noetherian AS regular \cite{AS}, then all results in the
paper apply to $A$.

To see so, we must find $K$ and $L$ as in Setup \ref{set:2} and show
that $\H(A)$ has a balanced dualizing complex.  The latter is true by
\cite[cor.\ 4.14]{Y}.

We know $\dim_k \Tor^{\H(A)}_*(k,k) < \infty$, and using the
Eilenberg-Moore spectral se\-quen\-ce shows
\[
  \dim_k \H(k \LTensor{A} k) < \infty
\]
whence $k$ is compact from either side.  Also, $\dim_k
\Tor^{\H(A)}_*(\H(A)^{\vee},k) = \dim_k \Ext_{\H(A)}^*(k,\H(A)) <
\infty$, and using the Eilenberg-Moore spectral sequence shows
\[
  \dim_k \H(\RHom_A(k,A)) = \dim_k \H(A^{\vee} \LTensor{A} k) < \infty
\]
whence $\sup \RHom_A(k,A) < \infty$ so $({}_{A}k)^* \in \langle k_A
\rangle$ by Lemma \ref{lem:bounded_above}.  Hence the $K$-projective
resolution of ${}_{A}k$ can be used for ${}_{A}K$.  Similarly, the
$k$-projective resolution of $k_A$ can be used for $L_A$.
\end{Example}

\begin{Example}
\label{exa:2}
If $A$ is commutative in the DG sense and $\H(A)$ is noetherian, then
all results in the paper apply to $A$.

To see so, we must again find $K$ and $L$ as in Setup \ref{set:2} and
show that $\H(A)$ has a balanced dualizing complex.  The former can be
done by using a DG version of the construction of the Koszul complex.

Since $\H(A)$ is graded commutative noetherian, it is a quotient of a
tensor product $B \underset{k}{\otimes} C$ where $B$ is a polynomial
algebra with finitely many generators in even degrees and $C$ is an
exterior algebra with finitely many generators in odd degrees.  It
follows from \cite[thm.\ 6.3]{VdB} that $B \underset{k}{\otimes} C$
has a balanced dualizing complex.  Now combine \cite[thm.\ 6.3]{VdB},
\cite[prop.\ 7.2(2)]{AZ}, and \cite[thm.\ 8.3(2+3)]{AZ} to see that so
does any quotient of $B \underset{k}{\otimes} C$.
\end{Example}

\begin{Example}
\label{exa:finite_dimensional}
If $\dim_k \H(A) < \infty$ then all results in the paper apply to $A$.
Namely, $\langle {}_{A}k \rangle = \langle {}_{A}A \rangle = \sD(A)$
and $\langle k_A \rangle = \langle A_A \rangle = \sD(A^{\opp})$ so we
can use ${}_{A}K = {}_{A}A$ and $L_A = A_A$.  Moreover, $\H(A)$ has
the balanced dualizing complex $\H(A)^{\vee}$.

In this case, we can easily find the dualizing DG module.  Since we
have ${}_{A}A \in \langle {}_{A}k \rangle$, the counit morphism
$\Gamma(A) \stackrel{\epsilon_A}{\longrightarrow} A$ is an
isomorphism.  Hence $F \LTensor{A} A \cong A$, that is, $F \cong A$,
and this is an isomorphism in $\sD(A^e)$; see the last paragraph of
Remark \ref{rmk:DG} and Remark \ref{rmk:two-sided}.  So the
dualizing DG module of $A$ is
\[
  D \cong A^{\vee}
\]
and the functors in Theorem \ref{thm:cor:two-sided} are just
$(-)^{\vee}$. 
\end{Example}

\begin{Example}
\label{exa:polynomials}
Let $A = k[T]$ have $T$ in cohomological degree $d \geq 1$ and
differential $\partial = 0$.  All results in the paper apply to $A$ by
Example \ref{exa:3}.  Let us compute the dualizing DG module.

While $A$ is commutative as a ring, it is not necessarily commutative
in the DG sense because this means $xy = (-1)^{|x||y|}yx$ for graded
elements $x$, $y$.  This fails if $d$ is odd and $k$ has
characteristic different from $2$.

However, it remains the case that $1 \mapsto T$ extends to a unique
homomorphism of DG $A$-bimodules $\Sigma^{-d}A
\stackrel{\varphi}{\longrightarrow} A$.  The homomorphism is injective
with cokernel ${}_{A}k_{A}$, so there is a distinguished triangle
\begin{equation}
\label{equ:tria1}
  \Sigma^{-d}A \stackrel{\varphi}{\longrightarrow} A \longrightarrow k
\end{equation}
in $\sD(A^e)$.

We can consider $N = k[T,T^{-1}]$ as a DG $A$-bimodule with $T$ in
cohomological degree $d$ and differential $\partial = 0$.  Then there
is a short exact sequence of DG $A$-bimodules $0 \rightarrow A
\rightarrow N \rightarrow C \rightarrow 0$ which induces a
distinguished triangle
\begin{equation}
\label{equ:tria2}
  \Sigma^{-1}C \rightarrow A \rightarrow N
\end{equation}
in $\sD(A^e)$.

It is easy to check that applying $\RHom_A(-,N)$ to \eqref{equ:tria1}
sends $\varphi$ to an isomorphism, so $\RHom_A(k,N) = 0$ whence
${}_{A}N \in \langle {}_{A}k \rangle^{\perp}$ and $\Gamma(N) = 0$; see
the last paragraph of Remark \ref{rmk:DG}.  We have
${}_{A}(\Sigma^{-1}C) \in \langle {}_{A}k \rangle$ by Lemma
\ref{lem:bounded_above}, so $\Gamma(\Sigma^{-1}C) \cong \Sigma^{-1}C$.
Note that this is an isomorphism in $\sD(A^e)$, see Remark
\ref{rmk:two-sided}, so it follows that applying $\Gamma$ to
\eqref{equ:tria2} produces an isomorphism $\Sigma^{-1}C \cong
\Gamma(A)$ in $\sD(A^e)$.  However, $\Gamma(A) = F \LTensor{A} A \cong
F$ in $\sD(A^e)$, so we get
\[
  \Sigma^{-1}C \cong F
\]
in $\sD(A^e)$.  Hence the dualizing DG module of $A$ is
\[
  D = F^{\vee} = (\Sigma^{-1}C)^{\vee}.
\]

More explicitly, $C$ is the DG quotient module $k[T,T^{-1}] / k[T]$,
and based on this, a concrete computation of $D =
(\Sigma^{-1}C)^{\vee}$ yields the following: As a graded vector space,
$D$ has a generator $e_\ell$ in cohomological degree $d\ell + d - 1$
for each $\ell \geq 0$.  It has differential $\partial = 0$, and the
left and right actions of $A$ on $D$ are given by
\[
  T^j e_\ell = e_{j+\ell}, \;\;\;
  e_\ell T^j = (-1)^{jd} e_{j+\ell}.
\]
As a DG left-$A$-module, $D$ is just $\Sigma^{-(d-1)}A$.  The right
action of $A$ is twisted by the DG algebra automorphism $\alpha : T^j
\mapsto (-1)^{jd}T^j$.  De\-no\-ting the twist by a superscript, we
finally have
\[
  D \cong (\Sigma^{-(d-1)}A)^{\alpha}
\]
in $\sD(A^e)$.

There is no way to get rid of the twist: If $D$ were isomorphic to
$\Sigma^{-(d-1)}A$ in $\sD(A^e)$, then the cohomologies of the two DG
modules would be isomorphic as graded $\H(A)$-bimodules.  However,
$\H(\Sigma^{-(d-1)}A)$ is a symmetric $\H(A)$-bimodule, but if $d$ is
odd and $k$ has characteristic different from $2$, then $\H(D)$ is
not.

The presence of twists in the theory of two sided duality is not
sur\-pri\-sing since it occurs already for rings, see for instance
\cite[thm.\ 7.18 and the remark preceding it]{Y}.
\end{Example}

\begin{Remark}
In Definition \ref{def:F}, the dualizing DG module $D$ depends on the
choice of the object $K$ made in Setup \ref{set:blanket}.  However, in
the two previous examples, the computations show that any choice of
$K$ produces the same $D$.  It would be interesting to know if $D$ is
unique in general.
\end{Remark}

\medskip
\noindent
{\bf Acknowledgement.}
Katsuhiko Kuribayashi provided useful comments on a preliminary
version, and Shoham Shamir kindly answered some questions on his paper
\cite{Shamir}.

\end{document}